\documentclass{amsart}
\usepackage{tipa}
\usepackage{amssymb}
\usepackage{stmaryrd}
\usepackage{amsmath,amsfonts,amsthm,mathrsfs,bbm,array,subfigure}

\newtheorem{theorem}{Theorem}[section]
\newtheorem{lemma}[theorem]{Lemma}

\newtheorem{corollary}[theorem]{Corollary}

\theoremstyle{definition}

\theoremstyle{remark}

\numberwithin{equation}{section}

\begin{document}

\title{A variety of Euler's conjecture}
\author{Tianxin Cai}
\address{Department of Mathematics, Zhejiang University,
Hangzhou, 310027, People's Republic of China }
\email{txcai@zju.edu.cn}
\thanks{Project supported by the Natural Science Foundation of China(\#10871169).}


\author{Yong Zhang}
\address{Department of Mathematics, Zhejiang University,
Hangzhou, 310027, People's Republic of China }

\email{zhangyongzju@163.com}


\date{}


\begin{abstract}
We consider a variety of Euler's conjecture, i.e., whether the
Diophantine system
\[\begin{cases}
n=a_{1}+a_{2}+\cdots+a_{s-1},\\
a_{1}a_{2}\cdots a_{s-1}(a_{1}+a_{2}+\cdots+a_{s-1})=b^{s}
\end{cases}\]
has solutions $n,b,a_i\in\mathbb{Z}^+,i=1,2,\ldots,s-1,s\geq 3.$ By
using the theory of elliptic curves, we prove that it has no
solutions $n,b,a_i\in\mathbb{Z}^+$ for $s=3$, but for $s=4$ it has
infinitely many solutions $n,b,a_i\in\mathbb{Z}^+$ and for $s\geq 5$
there are infinitely many polynomial solutions
$n,b,a_i\in\mathbb{Z}[t_1,t_2,\ldots,t_{s-3}]$ with positive value
satisfying this Diophantine system.
\end{abstract}

\maketitle
\section{Introduction}
In 1769, Euler conjectured that the Diophantine equation
\begin{equation}
a_1^s+a_2^s+\cdots+a_{s-1}^s=a_s^s
\end{equation}
has no positive integer solutions for $s\geq3$. It is called Euler's conjecture.

For $s=3$, (1.1) corresponds to the case $a_1^{3}+a_2^{3}=a_3^{3}$
of Fermat's Last Theorem and Fermat proved that it has no nontrivial
integer solutions.

For $s=4$, in 1988, N. Elkies [3] disproved Euler's conjecture by
showing that $a_1^4+a_2^4+a_3^4=a_4^4$ has
 infinitely many integer solutions. In particular, he gave the following solution:
\[2682440^{4}+15365639^{4}+18796760^{4}=20615673^{4},\]
and shortly after, R. Frye found the smallest counterexample:
\[95800^{4}+217519^{4}+414560^{4}=422481^{4}.\]

For $s=5$, in 1966, L. J. Lander and T. R. Parkin [5] found the
first counterexample:
\[27^{5}+84^{5}+110^{5}+133^{5}=144^{5}.\]
In 2004, J. Frye found the only other known primitive solution for
$s=5$:
\[55^{5}+3183^{5}+28969^{5}+85282^{5}=85359^{5}.\]

For $s\geq6$, there are no known solutions. More information about
Euler's conjecture can be found in [4]: \emph{D1 Sums of like
powers. Euler's conjecture.}

In 2011, the first author raised a new variant of the Hilbert-Waring Problem (cf. [1]): to express a positive integer $n$ as a sum of $s$ positive integers whose product is a $k$-th power, i.e.,
\begin{equation*}
  n=a_{1}+a_{2}+\cdots+a_{s}
\end{equation*}
such that
\begin{equation*}
  a_{1}a_{2}\cdots a_{s}=b^k
\end{equation*}
for $n,a_{i},b,k \in \mathbb{Z}^+.$

Now we expand this idea to Euler's conjecture and consider whether
the Diophantine system
\begin{equation}
\begin{cases}
n=a_{1}+a_{2}+\cdots+a_{s-1},\\
a_{1}a_{2}\cdots a_{s-1}(a_{1}+a_{2}+\cdots+a_{s-1})=b^{s}
\end{cases}
\end{equation}
has solutions $n,b,a_i\in\mathbb{Z}^+,i=1,2,\ldots,s-1,s\geq3.$

Obviously, the solutions of (1.1) is a subset of the solutions of (1.2). The motivation of studying (1.2) is that
we try to find a counterexample to Euler's conjecture for $s=6$. Although we can't get anyone, but we find some
interest results about this problem.

By using the theory of elliptic curves, we prove the following
theorems for (1.2).

\begin{theorem} For $s=3$, $(1.2)$ has no solution $n,b,a_i\in\mathbb{Z}^+,i=1,2$.
\end{theorem}

\begin{theorem} For $s=4$, $(1.2)$ has infinitely many solutions $n,b,a_i\in\mathbb{Z}^+,i=1,2,3$.
\end{theorem}

\begin{theorem} For $s\geq 5$, $(1.2)$ has infinitely many polynomial solutions
 $n,b,a_i\in\mathbb{Z}[t_1,t_2,\cdots,t_{s-3}]$
with positive value for $i\geq4$.
\end{theorem}

From Theorem 1.3, we have

\begin{corollary} For $s\geq 5$, $(1.2)$ has infinitely many solutions $n,b,a_i\in\mathbb{Z}^+,i=1,2,\cdots,s-1$.
\end{corollary}

\vskip10pt
\section{Two auxiliary lemmas}
To prove our theorems, we need the theory of the elliptic curves which includes Nagell-Lutz Theorem
 and the theorem of Poincar\'{e} and Hurwitz, we list them in the following two lemmas.

\begin{lemma}(Nagell-Lutz Theorem, see [7], p. 56) Let the equation of the elliptic curve be
 \[y^2=x^3+ax^2+bx+c~(a,b,c\in \mathbb{Z}),\]
the discriminant of the cubic polynomial is $\Delta=-4a^3c+a^2b^2+18abc-4b^3-27c^3,$
 let $P=(x,y)$ be a rational point of finite order, then $x$ and $y$ are integers;
  and either $y=0$ or else $y| \Delta$.
\end{lemma}

From this theorem, we know that if $x$ or $y$ is not an integer, then $P=(x,y)$ is a rational point
 of infinite order, hence there are infinitely many rational points on the elliptic curve.

\begin{lemma}(The theorem of Poincar\'{e} and Hurwitz, see [8], Chap. V, p.78, Satz 11) If the
 elliptic curve has infinitely many rational points, then it has infinitely many rational points
  in every neighborhood of any one of them.
\end{lemma}

\vskip10pt
\section{ Proofs of the Theorems}
In this section, we give the proofs of our theorems. For $s=3$,
there is no positive integer solution of (1.2), which is similar to
the Euler's conjecture. For $s=4$, there are infinitely many
positive integer solutions of (1.2).  But for $s\geq 5$, there are
infinitely many polynomial solutions
  of (1.2), which is different from the Euler's conjecture as we have
  known.\\

{\it \textbf{Proof of Theorem 1.1.}} For $s=3$, (1.2) is equal to
\[\begin{cases}
n=a_{1}+a_{2},\\
a_{1}a_{2}(a_{1}+a_{2})=b^{3},
\end{cases}\]
because of $n,b,a_i\in\mathbb{Z}^+,i=1,2$, from the second of the above equations, we have
\[\frac{a_{1}}{b}\frac{a_{2}}{b}\bigg(\frac{a_{1}}{b}+\frac{a_{2}}{b}\bigg)=1.\]
Let
\[b_i=\frac{a_{i}}{b}\in\mathbb{Q}^+,i=1,2,\]
we get
\[b_{1}b_{2}(b_{1}+b_{2})=1,\]
leading to
\[\bigg(\frac{b_{1}}{b_{2}}\bigg)^2+\frac{b_{1}}{b_{2}}=\frac{1}{b_{2}^3}.\]
Let \[u=\frac{b_{1}}{b_{2}},v=\frac{1}{b_{2}},\]we have
\[u^2+u=v^3,\]
by $y=16u+8,x=4v,$ we get
\[y^2=x^3+64,\]using the package of Magma [6], we can get the only rational points on it, i.e.,
\[(x,\pm y)=(8,24),(0,8),(-4,0).\]
Tracing back, we find that there is no integer solution of (1.2). \hfill $\oblong$\\

{\it \textbf{Proof of Theorem 1.2.}} For $s=4$, (1.2) is
\begin{equation}
\begin{cases}
n=a_{1}+a_{2}+a_{3},\\
a_{1}a_{2}a_{3}(a_{1}+a_{2}+a_{3})=b^{4},
\end{cases}
\end{equation}
because of $n,b,a_i\in\mathbb{Z}^+,i=1,2,3$, from the second equation of (3.1), we have
\[\frac{a_{1}}{b}\frac{a_{2}}{b}\frac{a_{3}}{b}\bigg(\frac{a_{1}}{b}+\frac{a_{2}}{b}+\frac{a_{3}}{b}\bigg)=1.\]
Let
\[b_i=\frac{a_{i}}{b}\in\mathbb{Q}^+,i=1,2,3,\]
we get
\[b_{1}b_{2}b_{3}(b_{1}+b_{2}+b_{3})=1.\]
It's easy to see that $(a_{1},a_{2},a_{3})=(1,2,24)$ satisfies (3.1), leading to
\[(b_{1},b_{2},b_{3})=\bigg(\frac{1}{6},\frac{1}{3},4\bigg).\]
Then
\begin{equation}
\begin{cases}
\begin{split}
&b_{1}b_{2}b_{3}=\frac{2}{9},\\
&b_{1}+b_{2}+b_{3}=\frac{9}{2}.
\end{split}
\end{cases}
\end{equation}

Next, we consider $b_i$ be unknowns, and will prove that there are infinitely
 many positive rational numbers satisfying
 (3.2).
Eliminating $b_3$ of $(3.2)$, we get
\[18b_{1}^2b_{2}+18b_{1}b_{2}^2-81b_{1}b_{2}+4=0,\]
leading to \[18\frac{b_{2}}{b_{1}}+18\bigg(\frac{b_{2}}{b_{1}}\bigg)^2-81\frac{b_{2}}{b_{1}}\frac{1}{b_{1}}+4\bigg(\frac{1}{b_{1}}\bigg)^3=0.\]
Taking \[u=\frac{b_{2}}{b_{1}},~v=\frac{1}{b_{1}},\]we have
\[18u^2+18u-81uv+4v^3=0.\]
Let \[y=384u-864v+192,~x=-32v+243,\]
we get
\[E:~~y^2=x^3-166779x+26215254.\]

By Lemma 2.1, to prove that there are infinitely many rational points on $E$, it is enough to find a rational point
 on $E$ with $x$-coordinate not in $\mathbb{Z}$. Using the package Magma [6], we can get the point
 $(30507/121,-584592/1331)$ with $x$-coordinate not in $\mathbb{Z}$, then there are infinitely many rational
 points on $E$.

From the above transformations, we have
\[
\begin{cases}
\begin{split}
&b_{1}=\frac{32}{243-x},\\
&b_{2}=\frac{-y+27x-6369}{12(243-x)},\\
&b_{3}=\frac{y+27x-6369}{243-x}.
\end{split}
\end{cases}
\]
Then
\[
\begin{cases}
\begin{split}
&a_{1}=\frac{32}{243-x}b,\\
&a_{2}=\frac{y-27x+6369}{12(243-x)}b,\\
&a_{3}=\frac{-y-27x+6369}{243-x}b
\end{split}
\end{cases}
\]
is a solution of (3.1).

To prove $b_i>0$, we get the equivalent condition
\[x<243,~|y|<27x-6369.\]
In virtue of Lemma 2.2, the elliptic curve has infinitely many rational points in every neighborhood of any one of them.
We should find a point satisfy the above equivalent condition. It's easy to see that the point $P=(235,8)$ satisfies it,
 then there are infinitely many rational points satisfying\[x<243,~|y|<27x-6369.\]
 Therefore, we can find infinitely many solutions in rational numbers $b_i>0,~i=1,2,3$ satisfying $(3.2),$ leading to
  integer numbers $a_i>0,~i=1,2,3,$ by multiply the least common denominator of $b_i>0,~i=1,2,3$. This proves that
  $(1.2)$ has infinitely many positive integer solutions for $s=4$. \hfill $\oblong$\\

Example, for $s=4$ the points
\[(x,y)=(235,8),\bigg(\frac{60266587}{257049},\frac{3852230624}{130323843}\bigg)\]
on the elliptic curve $y^2=x^3-166779x+26215254$, leading to
\[(a_1,a_2,a_3)=(1,2,24),(781943058,138991832,18609625).\]

{\it \textbf{Proof of Theorem 1.3.}} For $s\geq 5$, (1.2) is
\begin{equation}
\begin{cases}
n=a_{1}+a_{2}+\cdots+a_{s-1},\\
a_{1}a_{2}\cdots a_{s-1}(a_{1}+a_{2}+\cdots+a_{s-1})=b^{s},
\end{cases}
\end{equation}
from the second equation of (2.3), we have
\[\frac{a_{1}}{b}\frac{a_{2}}{b}\cdots \frac{a_{s-1}}{b}\bigg(\frac{a_{1}}{b}+\frac{a_{2}}{b}+\cdots+\frac{a_{s-1}}{b}\bigg)=1.\]
Let
\[b_i=\frac{a_{i}}{b}\in\mathbb{Q}^+,i=1,\cdots,s-1,\]
we get
\[b_{1}b_{2}\cdots b_{s-1}(b_{1}+b_{2}+\cdots+b_{s-1})=1.\]

For convenience, put $x=b_1,y=b_2,z=b_3$ and
\[u=b_{4}\cdots b_{s-1},v=b_{4}+\cdots+b_{s-1},\]
then we have\[xyzu(x+y+z+v)=1.\]

Let $z=ty$ in the above equation, we get
\begin{equation}
tuy^2x^2+ut(yt+y+v)y^2x-1=0,
\end{equation}
consider it as a quadratic equation of $x$, if it has rational
solutions, the discriminant
\[\Delta(y)=u^2t^2(t+1)^2y^4+2u^2v(t+1)t^2y^3+u^2v^2t^2y^2+4tu\]
should be a square. To prove Theorem 1.3, it is enough to show that
the set of $x\in \mathbb{Q}(t)$, such that (3.3) has a solution in
$\mathbb{Q}(t)$, is infinite. Then we need to show that there are
infinitely many $x\in \mathbb{Q}(t)$ such that the discriminant
$\Delta(y)$ should be a square in $\mathbb{Q}(t)$, which leads to
find infinitely many rational parametric solutions in
$\mathbb{Q}(t)$ on the following curve
\[C:~w^2=u^2t^2(t+1)^2y^4+2u^2v(t+1)t^2y^3+u^2v^2t^2y^2+4tu.\]

The discriminant of $C$ is
\[\Delta(t)=256(t+1)^4(64t^2+(128+v^4u)t+64)u^9t^9\] and is
non-zero as an element of $\mathbb{Q}(t)$ as $u,v\in \mathbb{Q}^+$.
Then $C$ is smooth.

By the Proposition 7.2.1 in [2], we can transform the curve $C$ into
a family of elliptic curves
\[
E:~Y^2=X(X^2+u^2v^2t^2X-16u^3t^3(t+1)^2),\] by the inverse
birational map $\phi:~(y,w)\longrightarrow(X,Y)$ with
\[y=\frac{Y-uvtX}{2ut(t+1)X},w=\frac{Y^2-u^2v^2t^2X^2-2X^3}{4ut(t+1)X^2},\] and
\[\begin{split}
X=&2ut(t+1)(ut(t+1)y^2+uvty-w),\\
Y=&2u^2t^2(t+1)^2(ut(t+1)y^2+uvty-w)(2(t+1)y+v).
\end{split}\]

To get the suitable points on $E$ such that we have infinitely many
rational solutions $y$, take $t=ut_0^2$ in the equation of $E$, we
get
\[E':~Y^2=X(X^2+u^4v^2t_0^4X-16u^6t_0^6(ut_0^2+1)^2).\]

Note that the point
\[P=(4u^3t_0^3(ut_0^2+1),4vu^5t_0^5(ut_0^2+1))\]lies on $E'$.
Using the group law on the elliptic curve, we obtain the point
\[[2]P=\bigg(\frac{16u^2t_0^2(ut_0^2+1)^2}{v^2},-\frac{64u^3t_0^3(ut_0^2+1)^3}{v^3}\bigg),\]
and the point $[4]P=p(t_0),q(t_0))$, where
\[\begin{split}p(t_0)=&\frac{u^2t_0^2(16u^2t_0^4+(32u+v^4u^2)t_0^2+16)^2}{64v^2(ut_0^2+1)^2},\\
q(t_0)=&-(u^3t_0^3(16u^2t_0^4+(32u+v^4u^2)t_0^2+16)(256u^4t_0^8\\
&+(1024u^3-64u^4v^4)t_0^6+(-u^4v^8-128u^3v^4+1536u^2)t_0^4\\
&+(-64v^4u^2+1024u)t_0^2+256))/(512v^3(ut_0^2+1)^3).\end{split}\]

Let us recall that on the elliptic curve
$y^2=x^3+a(t)x^2+b(t)x+c(t)$, where $a(t),b(t),c(t)\in
\mathbb{Z}[t]$, the points of finite order have coordinates in
$\mathbb{Z}[t]$. Therefore, to prove that the group
$E'(\mathbb{Q}(t_0))$ is infinite, it is enough to find a point with
coordinates not in $\mathbb{Z}[t]$.

Note that the $X$-coordinate of $[4]P$ is
\[\frac{u^2t_0^2(16u^2t_0^4+(32u+v^4u^2)t_0^2+16)^2}{64v^2(ut_0^2+1)^2},\]
when the numerator of the $X-$coordinate of it is divided by the
denominator, the remainder equals
\[r=u^3v^8(3ut_0^2+2).\]
For $u>0,v>0$, $r$ is not zero in $\mathbb{Q}(t_0)$, then the
$X-$coordinate of $[4]P$ is not a polynomial. Therefore, $[4]P$ is a
point of infinite order on $E'$. Hence, the group
$E'(\mathbb{Q}(t_0))$ is infinite.

Next, we will determine the positive polynomial solutions of (1.2).
From the birational map and the point $-[2]P$, i.e., the reflected
point of $[2]P$, we get
\[x=\frac{uv^3t_0}{2(4ut_0^2-uv^2t_0+4)},y=\frac{4ut_0^2-uv^2t_0+4}{2uvt_0(ut_0^2+1)},z=\frac{(4ut_0^2-uv^2t_0+4)t_0}{2v(ut_0^2+1)}.\]

To get $x>0,y>0,z>0$, from $u>0,v>0$, we need
\[
4ut_0^2-uv^2t_0+4>0,
\]
the discriminant of it is $\delta=u(uv^4-64).$ If $\delta<0,$ then
for any $t_0\in \mathbb{Q}$, we have
\[
4ut_0^2-uv^2t_0+4>0.
\]If $\delta>0,$ then
for any \[t_0\in
\bigg(0,\frac{uv^2-\sqrt{u^2v^4-64u}}{8u},\bigg),\bigg(\frac{uv^2+\sqrt{u^2v^4-64u}}{8u},\infty
\bigg),\] we have
\[
4ut_0^2-uv^2t_0+4>0.
\]
Therefore, for any $u>0,v>0,$ we can find infinitely many positive
rational numbers $t_0$ such that $x>0,y>0,z>0.$

Note that $x=b_1,y=b_2,z=b_3,$
\[u=b_{4}\cdots b_{s-1},v=b_{4}+\cdots+b_{s-1},\]and
\[b_i=\frac{a_{i}}{b}\in\mathbb{Q}^+,i=1,\ldots,s-1,\]
let \[t_1=t_0b',t_2=b_4b',t_3=b_5b',\ldots,t_{s-3}=b_{s-1}b',\]where
$b'$ is the least common  multiple of the denominator of
$b_i,i=1,\ldots,s-1$, then
\[a_i\in \mathbb{Z}[t_1,t_2,\ldots,t_{s-3}],i=1,\ldots,s-1\]and $a_i$ have positive
value, where $t_1=t_0b'$ is a positive parameter satisfying the
condition $4ut_0^2-uv^2t_0+4>0.$.

This completes the proof of Theorem 1.3. \hfill $\oblong$\\

Example, for $s=5$, we have
\[\begin{split}&x=\frac{uv^3t_0}{2(4ut_0^2-uv^2t_0+4)},y=\frac{4ut_0^2-uv^2t_0+4}{2uvt_0(ut_0^2+1)},\\
&z=\frac{(4ut_0^2-uv^2t_0+4)t_0}{2v(ut_0^2+1)},u=v=b_4,\end{split}\]
let $t_1=t_0,t_2=b_4,$ take
\[b=2uvt_0(ut_0^2+1)(4ut_0^2-uv^2t_0+4),\] then
\[\begin{split}
&a_1=t_1^2t_2^6(t_1^2t_2+1),a_2=(4t_1^2t_2-t_1t_2^3+4)^2,\\
&a_3=t_1^2t_2(4t_1^2t_2-t_1t_2^3+4)^2,a_4=2t_1t_2^3(t_1^2t_2+1)(4t_1^2t_2-t_1t_2^3+4).\end{split}\]
Let $t_1=1,t_2=1,$ we have
\[128=2+49+49+28,~2\cdot49\cdot49\cdot28\cdot(2+49+49+28)=28^5.\]
Let $t_1=2,t_2=1,$ we have
\[2000=20+324+1296+360,~20\cdot324\cdot1296\cdot360\cdot(20+324+1296+360)=360^5,\]which
can reduce to
\[500=5+81+324+90,~5\cdot81\cdot324\cdot90\cdot(5+81+324+90)=90^5.\]

\vskip10pt
\section{Some numerical solutions for $s=5,6$}
As we see in the above examples, the values of them are large, here
we list some smaller solutions in the following table for $s=5,6$ of
(1.2).

\begin{table}[htbp]
\caption{For $s=5,6$} \centering \subtable{
\begin{tabular}{cccccc}
\hline
$a_1$&$a_2$&$a_3$&$a_4$&$b$& $n$\\
\hline
$1$& $2$ & $12$  & $12$   & $6$ & $27$\\
\hline
$1$ & $4$  & $4$  & $18$ & $6$  & $27$\\
\hline
$1$ & $4$  & $20$  & $25$ & $10$  & $50$\\
\hline
$1$ & $3$  & $32$  & $36$ & $12$  & $72$\\
\hline
$1$ & $4$  & $12$  & $64$ & $12$  & $81$\\
\hline
$1$ & $3$  & $8$  & $96$ & $12$  & $108$\\
\hline
$1$ & $27$  & $36$  & $64$ & $24$  & $128$\\
\hline
$1$ & $1$  & $18$  & $108$ & $12$  & $128$\\
\hline
$1$ & $25$  & $54$  & $100$ & $30$  & $180$\\
\hline
$1$ & $4$  & $27$  & $256$ & $24$  & $288$\\
\hline
\end{tabular}
}
\qquad
\subtable{
\begin{tabular}{ccccccc}
\hline
$a_1$ & $a_2$  & $a_3$  & $a_4$ & $a_5$ & $b$ & $n$\\
\hline
$1$ & $1$  & $2$  & $2$ & $2$ & $2$ & $8$\\
\hline
$1$ & $6$  & $6$  & $6$ & $8$ & $6$ & $27$\\
\hline
$1$ & $1$  & $9$  & $9$ & $16$ & $6$ & $36$\\
\hline
$1$ & $2$  & $3$  & $12$ & $18$ & $6$ & $36$\\
\hline
$1$ & $9$  & $12$  & $18$ & $24$ & $12$ & $64$\\
\hline
$1$ & $4$  & $16$  & $24$ & $27$ & $12$ & $72$\\
\hline
$1$ & $6$  & $9$  & $24$ & $32$ & $12$ & $72$\\
\hline
$1$ & $4$  & $8$  & $32$ & $36$ & $12$ & $81$\\
\hline
$1$ & $4$  & $12$  & $16$ & $48$ & $12$ & $81$\\
\hline
$1$ & $2$  & $9$  & $36$ & $48$ & $12$ & $96$\\
\hline
\end{tabular}
}
\end{table}

\vskip20pt
\bibliographystyle{amsplain}

\end{document}